\documentclass[a4paper,12pt]{article}
\usepackage{times, url}
\usepackage{amsmath,amssymb,amsthm,amsfonts,bm,float,cite}
\usepackage{booktabs}
\usepackage{boldline}

\usepackage{hyphenat}
\PassOptionsToPackage{hyphens}{url}
\usepackage{hyperref}

\usepackage[center]{caption}

\newtheorem{definition}{Definition}[section]
\newtheorem{theorem}[definition]{Theorem}

\newtheorem{manualconjectureinner}{Conjecture}
\newenvironment{manualconjecture}[1]{%
 \IfBlankTF{#1}
  {\renewcommand{\themanualconjectureinner}{\unskip}}
  {\renewcommand\themanualconjectureinner{#1}}%
 \manualconjectureinner
}{\endmanualconjectureinner}

\usepackage{mathtools}
\DeclarePairedDelimiter\ceil{\lceil}{\rceil}
\DeclarePairedDelimiter\floor{\lfloor}{\rfloor}

\makeatletter
\renewcommand{\@biblabel}[1]{[#1]\hfill}
\makeatother

\begin{document}
\setcounter{page}{1}

\begin{center}
{\LARGE \bf Minimum-Link Covering Trails for any Hypercubic Lattice\\ } 
\vspace{10mm}

{\Large Marco Rip\`a}
\vspace{2mm}

World Intelligence Network \\ 
Rome, Italy \\
e-mail: \url{marcokrt1984@yahoo.it}
\vspace{8mm}

\end{center}


\noindent {\bf Abstract.} In 1994, Kranakis et al. published a conjecture about the minimum link-length of every rectilinear covering path for the $k$-dimensional grid $P(n,k) \coloneqq \{0,1, \dots, n-1\} \times \{0,1, \dots, n-1\} \times $ \linebreak$ \cdots \times \{0,1, \dots, n-1\}$. In this paper, we consider the general, NP-complete, Line-Cover problem, where the edges are not required to be axis-parallel, showing that the original Theorem 1 by Kranakis et al. no longer holds when the aforementioned constraint is disregarded. Furthermore, for any $n$ greater than two, as $k$ approaches infinity, the link-length of any minimal (non-rectilinear) polygonal chain does not exceed Kranakis' conjectured value of $\frac{k}{k-1} \cdot n^{k-1}+O(n^{k-2})$ only if we introduce a multiplicative constant $c \geq 1.5$ for the lower order terms (e.g., if we select $n=3$ and assume that $c<1.5$, starting from a sufficiently large $k$, it is not possible to visit all the nodes of $P(n,k)$ with a trail of link-length $\frac{k}{k-1} \cdot n^{k-1}+c \cdot n^{k-2}$).

\vspace{8mm}


\section{Introduction} \label{sec:Intr}
\sloppy We take into account a multidimensional extension of a well-known puzzle involving lateral thinking \cite{Kihn:4, Loyd:7, Lung:8} and examine a generalization of an intriguing conjecture by Kranakis, Krizanc, and Meertens, published in 1994 on Ars Combinatoria, vol. 38, p. 191 \cite{Kranakis:6}.

In order to compactly describe the related open problems and the main results of the present paper, let us give a few basic definitions first (we will denote $\mathbb{N}$ as the set of positive integers $\{1,2,3,\ldots\}$).

\begin{definition} \label{def1.1}
Let $n, k \in \mathbb{N}$. Accordingly, $P(n,k) \coloneqq \{0,1, \dots, n-1\} \times \{0,1, \dots, n-1\} \times \cdots \times \{0,1, \dots, n-1\}$ is a finite set of $n^k$ points in $\mathbb{R}^{k}$.
\end{definition}

\begin{definition} \label{def1.2}
We define as $h(n,k)$ the number of edges of the minimum-link covering trail for $P(n,k)$ (i.e., the link-length of any polygonal chain that visits each node of $P(n,k)$ at least once and that is characterized by the minimum number of edges). Furthermore, let us call $h_l(n,k)$ and $h_u(n,k)$ the proven lower and upper bound for $h(n,k)$ (respectively), so $h_l(n,k) \leq h(n,k) \leq h_u(n,k)$ holds for every $P(n,k) \subset \mathbb{R}^{k}$.
\end{definition}

\begin{definition} \label{def1.3}
Let $k,n_1,n_2, \dots, n_k \in \mathbb{N}$ be such that $n_1 \leq n_2 \leq \dots \leq n_k$. Accordingly, $G(n_1, n_2, \dots, n_k) \coloneqq \{0,1, \dots, n_1-1\} \times \{0,1, \dots, n_2-1\} \times \dots \times \{0,1, \dots, n_k-1\}$ is a \linebreak finite set of $\prod_{i=1}^{k} n_i$ points in $\mathbb{R}^{k}$.
\end{definition}

\begin{definition} \label{def1.4}
We define as $h(n_1, n_2, \dots, n_k)$ the number of edges of the minimum-link covering trail for $G(n_1, n_2, \dots, n_k)$. Furthermore, let us call $h_l(n_1, n_2, \dots, n_k)$ a proven lower bound for $h(n_1,n_2, \dots, n_k)$ and, similarly, let us indicate as $h_u(n_1, n_2, \dots, n_k)$ a valid upper bound for the same set of $n_1 \cdot n_2 \cdots n_k$ points.
\end{definition}

The NP-complete problem (see \cite{Bereg:1}, p. 514) that we aim to solve asks to visit, at least once, all the nodes of $P(n,k)$ with a polygonal chain that is characterized by the minimum number of edges \cite{Moret:3, Kranakis:6}.

In 1994, Evangelos Kranakis et al. \cite{Kranakis:6} conjectured that $h(n,k) \leq \frac{k}{k-1} \cdot n^{k-1}+O(n^{k-2})$ for each $k \in \mathbb{N}-\{1,2\}$, under two additional conditions: all line-segments are axis-parallel \cite{Collins:2}, and every point of the $k$-dimensional grid should not be visited twice.

Kranakis' expected value of the link-length of the minimal rectilinear paths could not be smaller than any proven lower bound for the same optimization problem, and this proposition is automatically confirmed if the aforementioned additional constraints are removed. For example, in the unconstrained case, it has been constructively shown that $h(4,3) \leq 23 < \frac{k}{k-1} \cdot n^{k-1}$ \cite{Ripa:11} and $h(5,3) \leq 36 < \frac{k}{k-1} \cdot n^{k-1}$ \cite{Ripa:9}. Thus, we will demonstrate that, for any $k \geq 8$ and $n=3$, $h(n,k)$ lies between Kranakis' conjectured bound of $\frac{k}{k-1} \cdot n^{k-1} + n^{k-2}$ and Bereg's proven upper bound of $\frac{k}{k-1} \cdot n^{k-1} + n^{k-\frac{3}{2}}$ \cite{Bereg:1}.


\section{Current upper bound vs. best theoretical solution}

Let $k,n_1 \in \mathbb{N}-\{1,2\}$ and assume $n_1 = n_k = n$. From \cite[Equation (9)]{Ripa:10}, and given that $\frac{n_k+n_{k-1}}{2} = n$ by hypothesis, it follows that the trivial lower bound
\begin{equation} \label{eq:1}
h_l(n,k) = \ceil*{\frac{n^k-(k-2) \cdot n^2+(k-2) \cdot n-n+n \cdot ((k-2) \cdot n-k+2)}{n-1}}+1=\frac{n^k-1}{n-1}.
\end{equation}

Generally speaking, as shown by Theorem \ref{Theorem 2.1} below, if $3 \leq n_1 \leq n_2 \leq \cdots \leq n_k$ and $k \geq 2$, then it is not possible to visit all the nodes of $G(n_1,n_2, \dots, n_k)$ (see Definition \ref{def1.3}) with a covering trail whose link-length is smaller than $\ceil*{2 \cdot \frac{\prod_{i=1}^{k} n_i - n_k}{n_k + n_{k-1} - 2}}+1$.

\begin{theorem} \label{Theorem 2.1}
If $k > 1$ and $n_1 > 2$, then
\begin{equation} \label{eq:2}
h(n_1, n_2, \dots, n_k) \geq \ceil*{2 \cdot \frac{\prod_{i=1}^{k} n_i - n_k}{n_k + n_{k-1} - 2}} + 1
\end{equation}
holds for each $k$-tuple $(n_1, n_2, \ldots, n_k)$.
\end{theorem}

\begin{proof} Let the $k$ integers $n_1$, $n_2$, $\dots$, $n_k$ be such that $3 \leq n_1 \leq n_2 \leq \cdots \leq n_k$ (see Definition \ref{def1.3}). 

We immediately notice that it is not usually possible to intersect more than $(n_k-1)+(n_{k-1}-1)= n_k+n_{k-1}-2$ points using two (consecutive) straight lines connected at their endpoints; however, there is one exception (which, for the sake of simplicity, we may assume it as in the case of the first line drawn). In this circumstance, it is possible to fit $n_k$ points with the first line segment, then $n_{k-1}-1$ points using the second line and $n_k-1$ points with the third one (or, alternatively, $n_{k-1}-2$ points using the second line and $n_k$ points with the next one), and so forth.

In general, since for any $t \in \mathbb{N}$ we have $\frac{(n_{k-1}-1)+(n_k-1)}{2} \geq \frac{{\ceil*{\frac{t}{2}}} \cdot (n_{k-1}-1) + \floor*{\frac{t}{2}} \cdot (n_k-1)}{t}$, there exists no polygonal chain of $1+t$ links that visits more than $n_k + \frac{(n_k-1)+(n_{k-1}-1)}{2} \cdot t$ nodes of $G(n_1, n_2, \ldots, n_k)$.

Thus, $h(n_1, n_2, \dots, n_k)$ must satisfy the relation
\begin{equation} \label{eq:3}
\prod_{i=1}^{k} n_i \leq n_k + \frac{n_k + n_{k-1} - 2}{2} \cdot \left(h(n_1, n_2, \dots, n_k) - 1 \right).
\end{equation}

Hence,
\begin{equation*}
h(n_1, n_2, \dots, n_k) - 1 \geq 2 \cdot \frac{\prod_{i=1}^{k} n_i - n_k}{n_k + n_{k-1}-2}.
\end{equation*}

Since $h(n_1, n_2, \dots, n_k)$ is a natural number by definition, we can finally set
\begin{equation} \label{eq:4}
h_l(n_1, n_2, \dots, n_k) \coloneqq \ceil*{2 \cdot \frac{\prod_{i=1}^{k} n_i - n_k}{n_k + n_{k-1} - 2}} + 1
\end{equation}
so that $h(n_1, n_2, \dots, n_k) \geq h_l(n_1, n_2, \dots, n_k)$,
and this concludes the proof of Theorem \ref{Theorem 2.1}.	
\end{proof}

Now, let $n \coloneqq n_1= n_2= \cdots =n_k$ belong to the set $\mathbb{N}-\{1, 2\}$, as usual.

Since $h(n,3) \leq \floor*{\frac{3}{2} \cdot n^2} - \floor*{\frac{n-1}{4}} + \floor*{\frac{n+1}{4}} - \floor*{\frac{n+2}{4}} + \floor*{\frac{n}{4}} + n - 2$ by \cite{Bencini:12}, it follows that
\vspace{-1.5mm}
\begin{equation} \label{eq:5}
h(n,k) \leq h_u(n,k) \coloneqq \Bigg( \floor*{\frac{3}{2} \cdot n^2} - \floor*{\frac{n-1}{4}} + \floor*{\frac{n+1}{4}} - \floor*{\frac{n+2}{4}} + \floor*{\frac{n}{4}} + n - 1 \Bigg) \cdot n^{k-3} - 1
\end{equation}
holds for every $k > 2$. 

Thus,
\vspace{-1mm}
\begin{equation} \label{eq:6}
\begin{gathered}
\lim_{k\to\infty} \frac{n^k}{h_u(n,k)} = \lim_{k\to\infty} \frac{n^k}{\Big( \floor*{\frac{3}{2} \cdot n^2} - \floor*{\frac{n-1}{4}} + \floor*{\frac{n+1}{4}} - \floor*{\frac{n+2}{4}} + \floor*{\frac{n}{4}} + n - 1 \Big) \cdot n^{k-3} - 1} \\ \geq \lim_{k\to\infty} \frac{n^k}{\big( \frac{3 \cdot n^2}{2} - \frac{n-1}{4} + \frac{n+1}{4} - \frac{n+2}{4} + \frac{n}{4} + n \big) \cdot n^{k-3}} =
\lim_{k\to\infty} \frac{2 \cdot n^k}{n^{k-3} \cdot (3 \cdot n^2 + 2 \cdot n)} = \frac{2 \cdot n^2}{3 \cdot n + 2},
\end{gathered}
\end{equation}

\noindent whereas
\begin{equation} \label{eq:7}
\lim_{k\to\infty} \frac{n^k}{h_l(n,k)} = \lim_{k\to\infty} \frac{n^k}{\frac{2 \cdot n^k - n + n - 2}{n + n - 2}} = \lim_{k\to\infty} \frac{n^k}{\frac{n^k-1}{n-1} } = n-1.
\end{equation}

\vspace{2mm}
On the other hand, it is easy to check that, for any $n > 2$,

\centerline{ $\lim\limits_{k\to \infty} \frac{n^k}{h_u(n,k)} = \lim\limits_{k\to \infty} \frac{n^k}{\big( \floor*{\frac{3 \cdot n^2}{2}} - \floor*{\frac{n-1}{4}} + \floor*{\frac{n+1}{4}} - \floor*{\frac{n+2}{4}} + \floor*{\frac{n}{4}} + n - 1 \big) \cdot n^{k-3} - 1} \Rightarrow \lim\limits_{k\to \infty} \frac{n^k}{h_u(n,k)} < n-1$.}

\vspace{2mm}
It follows that, on average, as $k$ approaches infinity, the efficiency loss per link is equal to
\begin{equation*}
L(n) \coloneqq \lim\limits_{k\to \infty} \frac{n^k}{h_l(n,k)} - \lim\limits_{k\to \infty} \frac{n^k}{h_u(n,k)}
\end{equation*}
unvisited points.

Hence,
\begin{equation} \label{eq:8}
L(n) \leq \lim_{k\to\infty} \bigg( n - 1 + \frac{1}{h_l(n,k)} \bigg) - \frac{2 \cdot n^2}{3 \cdot n + 2} = \frac{n^2 - n - 2}{3 \cdot n + 2}.
\end{equation}

This is reasonable because, by Theorem \ref{Theorem 2.1}, for any $n$ greater than $2$, we know that $n-1+\frac{1}{h_l(n,k)}$ is the maximum average number of ``new'' nodes of $P(n,k)$ visited by each edge of a trail since the highest theoretical number of unvisited nodes covered by $t$ consecutive edges of a trail is $t \cdot (n-1)+1$ (for any such $n$).

In particular, if $n=3$, we can improve (\ref{eq:1}) as (\ref{eq:9}) (see \cite{Ripa:11})

\begin{equation} \label{eq:9}
h_l(3,k)=h(3,k)=\frac{3^k-1}{2}.
\end{equation}

Now, we observe that
\begin{equation} \label{eq:10}
h(3,k)<3^{k-1}+\frac{3}{2} \cdot 3^{k-2}	
\end{equation}
holds for each positive integer $k$, and Kranakis' conjectured upper bound $\frac{k}{k-1} \cdot n^{k-1}+O(n^{k-2})$ (valid for $n \geq 3$ and $k \geq 4$), concerning any minimal rectilinear walk covering $P(n,k)$, implies the existence of a constant $c \geq \frac{3}{2}$ such that $\frac{k}{k-1} \cdot n^{k-1}+c \cdot n^{k-2} \geq h(n,k)$.

Moreover, we note that
\begin{equation} \label{eq:11}
\lim_{k\to\infty} \frac{3^k}{h(3,k)} = \lim_{k\to\infty} \frac{3^k}{\frac{3^k-1}{2}} = 2,
\end{equation}
indicating zero efficiency loss per link since $L(3)=\lim\limits_{k\to\infty} \frac{3^k}{h_l(3,k)} - \lim\limits_{k\to\infty} \frac{3^k}{h_u(3,k)}=2-2 = 0$, as opposed to the generic bound $L(3) \leq \frac{4}{11}$ derived from (\ref{eq:8}).

\begin{manualconjecture}{2.1} \label{Conjecture 2.1}
For any given $P(n,k)$, we conjecture that the edges of each minimal covering trail visit (on average) less than $n-1$ new nodes if and only if $n \in \mathbb{N}-\{1, 2, 3\}$ and $k \in \mathbb{N}-\{1\}$.
\end{manualconjecture}

\noindent \textbf{Remark 2.1} \label{Remark 2}
If Conjecture \ref{Conjecture 2.1} holds, then there exists no $n \in \mathbb{N}-\{1,2,3\}$ such that \allowbreak$h(n,k)< \frac{n^k}{n-1}$, and hence $\frac{n^k}{h(n,k)} \geq n-1$ implies $n \leq 3$. This is certainly true for $(n,k)=(2,3)$, since $(1-\sqrt{2}, -\sqrt{2},0) \rightarrow (\sqrt{2},\sqrt{2},0) \rightarrow (\frac{1}{2},\frac{1}{2},2 \cdot \sqrt{3}-\sqrt{\frac{3}{2}}) \rightarrow (\sqrt{2},1-\sqrt{2},0) \rightarrow (1-\sqrt{2},\sqrt{2},0) \rightarrow \linebreak (\frac{1}{2},\frac{1}{2},2 \cdot \sqrt{3}-\sqrt{\frac{3}{2}}) \rightarrow (1-\sqrt{2},1-\sqrt{2},0)$ is a covering cycle for $P(2,3)$ (see \cite{Ripa:13}, pp. 162--163) and, by definition, it follows that $h_u(2,3)=6 \Rightarrow \frac{2^3}{h_u(2,3)} \leq \frac{2^3}{h(2,3)} \Rightarrow \frac{8}{h(2,3)}>2-1$. Likewise, (\ref{eq:9}) underlines that Kranakis' Theorem 1 of \cite{Kranakis:6} is not valid if we disregard the rectilinear walk constraint and allow covering paths as the above.


\section{Conclusion}
We have shown that Kranakis' conjectured upper bound of $\frac{k}{k-1} \cdot n^{k-1}+O(n^{k-2})$, for minimal rectilinear walks covering $P(n,k)$, can be rewritten, taking into account the general covering trails considered in \cite{Ripa:11, Bencini:12}, as $\frac{k}{k-1} \cdot n^{k-1} + c \cdot n^{k-2}$, where $c \geq \frac{3}{2}$. On the contrary, the upper bound proved by Bereg et al. \cite{Bereg:1} of $\frac{k}{k-1} \cdot 3^{k-1}+3^{k - \frac{3}{2}}$ definitely holds for every finite set of $3^k$ points in $\mathbb{R}^{k}$.

Lastly, for any nontrivial value of $k$, we see that $\frac{2^k}{h(2,k)} > \frac{3^k}{2 \cdot h(3,k)} > \frac{4^k}{3 \cdot h(4,k)}$ and Conjecture \ref{Conjecture 2.1} leads to a more general research question, which can be formulated as follows:

\centerline{``Does the relation $\frac{n^k}{(n-1) \cdot h(n,k)} > \frac{(n+1)^k}{n \cdot h(n+1,k)}$ hold for any $n,k \in \mathbb{N}-\{1\}$?''}

\noindent (e.g., the answer is affirmative if we select $k=2$ since $n>2 \Rightarrow h(n,2)=2 \cdot n-2$ by \cite{Keszegh:5}, and so $\frac{n^2}{(n-1) \cdot h(n,2)}= \frac{1}{2 \cdot (n-1)^2} + \frac{1}{n-1} + \frac{1}{2}$ proves that $\frac{n^2}{(n-1) \cdot h(n,2)} > \frac{(n+1)^2}{n \cdot h(n+1,2)}$ holds for every $n$).

\bibliographystyle{plain}
\bibliography{Minimum_link_covering_trails_for_any_hypercubic_lattice}

\end{document}